\documentclass{article}

\usepackage{amssymb}

\begin{document}
\input amssym.def

\renewcommand\subsection[1]
{\refstepcounter{subsection}
{\nopagebreak\smallskip\nopagebreak
\bf \arabic{section}.\arabic{subsection}.
{#1}}}

\renewcommand{\thesubsection}{\arabic{section}.\arabic{subsection}}

\renewcommand{\theequation}{\arabic{section}.\arabic{equation}}

\def\B{{\cal{B}}}
\def\C{{\Bbb{C}}}
\def\d{{\partial}}
\def\heis{{\frak{h}}}
\def\K{{\cal{K}}}
\def\L{{\cal{L}}}
\def\R{{\Bbb{R}}}
\def\S{{\cal S}}
\def\T{{\Bbb{T}}}
\def\tK{{\tilde{\K}}}
\def\tQ{{\tilde{Q}}}
\def\tu{{\tilde{u}}}
\def\tv{{\tilde{v}}}
\def\Z{{\Bbb{Z}}}

\newtheorem{theorem}{Theorem}[section]
\newtheorem{prop}[theorem]{Proposition}
\newtheorem{lemma}[theorem]{Lemma}


\def\SP{\mathbf{Sp}}
\def\Sp{\mathrm{Sp}}
\def\graph{\mathrm {graph}}
\def\Mor{\mathrm{Mor}}
\def\cS{{\cal S}}
\def\cK{{\cal K}}
\def\cM{{\cal M}}
\def\fZ{{\frak Z}}
\def\End{\mathrm{End}}
\def\cL{{\cal L}}
\def\cD{{\cal D}}
\def\cQ{{\cal Q}}
\def\cB{{\cal B}}

\newfont{\ya}{cmff10}
\def\add{\ttfamily}

\def\rr{\rightrightarrows}

\def\ov{\overline}

\renewcommand{\Re}{\mathop{\mathrm{Re}}\nolimits}
\newcommand{\rk}{\mathop{\mathrm{rk}}\nolimits}

\begin{center}

{\bf \large

Zak transform, Weil representation,
and integral  operators with theta-kernels}

\medskip

Foth T.\footnote{Department of Mathematics, University
of Michigan, Ann Arbor, MI 48109, USA

\tt foth@umich.edu},
Neretin Yu.A.\footnote{ITEP, B.Cheremushkinskaya 25,
Moscow 117 259, Russia

\tt neretin@mccme.ru}

\end{center}

\medskip

{ \small
{\sc Abstract.} The Weil representation of a real
symplectic group $\Sp(2n,\R)$ admits
a canonical extension to
a holomorphic representation of a certain
complex semigroup consisting of Lagrangian
linear relations (this semigroup includes the
Olshanski
semigroup). We obtain the explicit realization
of the Weil representation of this semigroup
in the Cartier model, i.e., in the space of smooth
sections of a certain line bundle on the
$2n$-dimensional
torus $\T^{2n}$. We show that operators of
the representation are integral operators
whose kernels are theta-functions on $\T^{4n}$.

We also extend this construction to a functor
from a certain category of Lagrangian linear relations
between symplectic vector spaces of different
dimensions
to a category of integral operators
acting on sections of line bundles on the tori.
 }

\medskip

\section{Introduction}

\subsection{Theta-kernels.}
Let $T$ be an $\alpha\times \alpha$ symmetric
matrix with
negative definite real part. Let $z$  range in
$\R^\alpha$.
We define a {\it theta-function} $\theta(T;z;\zeta)$
 as a function
in variables $(z,\zeta)\in\R^\alpha \times \R^\alpha$
given by
\begin{equation}
\theta[T;z;\zeta]=
\sum_{k\in \Z^\alpha} \exp\bigl\{\frac 12 (z+ 2\pi
k)^t T
(z+2\pi k)  + ik\cdot\zeta
                   \bigr\}
\end{equation}

This expession  is a variant of a multivariate
theta-function, see \cite{Mum1}, Chapter 2, Section 1
(see definition of theta-functions with characteristics).

Consider the space  $L^2([0,2\pi]^{2n})$.
We write its elements as functions $f(x,\xi)$,
where $x$ and $\xi$ are elements of $[0,2\pi]^n$.
We consider integral operators in
$L^2([0,2\pi]^{2n})$, with $\theta$-kernels.
More precisely, let
 $\alpha=2n$, $z=\pmatrix{x&y}^t\in\R^n\oplus\R^n$,
$\zeta=(\xi,-\eta)^t\in \R^n\oplus \R^n$.
Let $T=\pmatrix{A&B\cr B^t&C}$ be a
$(n+n)\times (n+n)$ matrix.
We consider integral   operators in
$L^2([0,2\pi]^{2n})$
 having the form
\begin{equation}
A f(x,\xi)=\frac{1}{(2\pi)^{3n/2}}
\int_{[0,2\pi]^n \times [0,2\pi]^n}
\theta\left[
\pmatrix{ A&B\cr B^t&C};
x,y; \xi, \eta
\right]\, f(y,\eta)\, dy\,d\eta
\label{theform}
\end{equation}
We show that  these operators form a semigroup,
 i.e., the  product of two such integral operators
has the same
form (see Proposition \ref{dubl}).
In fact, multiplying two integral operators associated
with
matrices $\pmatrix{ A&B\cr B^t&C}$ and
$\pmatrix{ K&L\cr L^t&M}$ we obtain the
integral operator associated
with the matrix
\begin{equation}
\pmatrix{A-B(C+K)^{-1}B^t & -B(K+C)^{-1}L \cr
-L^t(C+K)^{-1}B^t & M-L^t(C+K)^{-1}L}.
\label{duble}
\end{equation}


\subsection{Lagrangian linear relations.}
The formula (\ref{duble}) for a strange
multiplication of matrices
admits a transparent algebraic interpretation.
In fact, the semigroup of matrices with
the multiplication (\ref{duble})
 is isomorphic to a certain semigroup of
Lagrangian linear relations (see its
description below, in Section \ref{prelim}).

 This semigroup of Lagrangian linear
relations   has several interesting realizations by
integral
operators.

First, it  is isomorphic to the semigroup of integral
operators in $\R^{n}$ having
 the form
\begin{equation}
\B f(x)=\frac 1{(2\pi)^{n/2}}
\int_{\R^n}
\exp\Bigl\{ \frac 12 \pmatrix{x^t & y^t}
                    \pmatrix{A&B\cr B^t&C}
                    \pmatrix{x\cr y}\Bigr\}
\,f(y)\,dy
\end{equation}

Equivalently,  this  semigroup is isomorphic
to the semigroup of Gauss
operators
in the holomorphic (Segal--Bargmann)
model of the boson Fock space with $n$ degrees of
freedom
(see \cite{NNO}, \cite{Ner1}, later this was obtained
in \cite{Hor}).

An infinite-dimensional variant of this semigroup
 has a realization in the space of symmetric
functions, see \cite{Ner4}.

All these realizations are different models of the
Weil representation.%
\footnote{Other terms are the {\it
Segal--Shale--Weil
representation,
the harmonic representation, the oscillator
representation.}
Apparently, it was discovered by K.O.Friedrichs
(see \cite{Fri}) near
1950. He formulates a correct theorem on the
infinite
dimensional
symplectic group and proposes a proof "at the physical
level";
it is not satisfactory from the mathematical point
of
view.
But for finite dimensional symplectic groups
his arguments are sufficient.}
A relation of the Weil representation
with theta-functions seems well known, see
\cite{Car}, \cite{LV}, \cite{Mum3}.
A possibility to realize the Weil representation
by integral operators with theta-kernels was
conjectured
by R.Howe (private discussion, 1994). Authors are
grateful
to R.Howe, A.A.Kirillov, and E. Kaniuth
for discussions of the subject.

\subsection{Structure of this paper.}
Our main result (Theorem  \ref{main2})
is the multiplication law
for the operators (\ref{theform}) with
$\theta$-kernels.
We also show
(Theorem \ref{main1}) that the Zak transform
maps Gauss operators in $L^2(\R^n)$ to operators
with $\theta$-kernels in $L^2(\T^{2n})$.
 These theorems are  proved in Section
\ref{opertheta}.

Section \ref{prelim} contains preliminaries on Gauss
operators
(see also \cite{How}, \cite{Fol},
\cite{Ner2}-\cite{Ner3})
and Lagrangian linear relations
(see also  \cite{Ner2}).

\section{ Preliminaries.
Gauss operators in $L^2(\R^n)$
and the category of Lagrangian linear relations.}
\label{prelim}

\setcounter{equation}{0}

\subsection{Category of Gauss operators.}
Denote by $\cD(\R^n)$ the space of
complex-valued
$C^\infty$-functions
on $\R^n$ with compact support.
Denote by $\cS(\R^n)$ the Schwartz space of smooth
rapidly decreasing complex-valued functions on $\R^n$,
denote by $\cS'(\R^n)$ and $\cD'(\R^n)$
 their dual spaces, i.e., the space
of all distributions on $\R^n$ and the space
of
tempered distributions.

Fix integers $m,n\ge 0$.
Denote by $\Omega_{m,n}$ the set of
$(n+m)\times (n+m)$ block matrices of the form
$S=\pmatrix{K & L \cr L^t & M}$, where

\smallskip

1) $S=S^t$, i.e., $K^t=K$, $M^t=M$,

 2) $(-\Re S)$ is strictly positive definite.

\smallskip

For a matrix $S\in\Omega_{m,n}$,
define the integral
operator $L^2(\R^m)\to L^2(\R^{n})$ of the form
\begin{equation}
\cB[S]f(x)=\frac 1{(2\pi)^{(m+n)/4}}
\int_{\R^m}\exp\Bigl\{\frac{1}{2}\pmatrix{x^t &
y^t}S
\pmatrix{x \cr y}\Bigr\}f(y)dy.
\label{intop}
\end{equation}

Since the kernel of the operator is square integrable,
$\cB[S]$ is a Hilbert--Schmidt operator $L^2(\R^m)\to
L^2(\R^n)$.
Since the kernel is an element of $\cS(\R^{n+m})$,
our operator is a bounded operator
$\cS'(\R^m)\to\cS(\R^n)$ (see
a variant of the Kernel Theorem
 \cite{Hor0}, Theorem 5.2.6).

It can be easily checked by a direct calculation
(see Howe \cite{How}%
\footnote{To be precise, he considers the case
$m=n$})
that the product of two operators of this form
has the same form, the precise statement is as
follows.

\begin{theorem}  \label{Th-gauss-prod}
Let
$$S_1=\pmatrix{K_1 & L_1 \cr L_1^t &
M_1}\in\Omega_{n,k}
,\qquad
S_2=\pmatrix{K_2 & L_2 \cr L_2^t & M_2}\in\Omega_{m,n}
$$
then
the composition is given by
\begin{equation}
\cB[S_1] \cB[S_2]=\lambda(S_1,S_2)\cB[S_3]
\label{kompop}
\end{equation}
where
\begin{equation}
\lambda(S_1,S_2)=\det [(-M_1-K_2)^{-1/2}]
\label{lamb}
\end{equation}
and
\begin{equation}
S_3=\pmatrix{K_1-L_1(M_1+K_2)^{-1}L_1^t &
-L_1(M_1+K_2)^{-1}L_2 \cr
-L_2^t(M_1+K_2)^{-1}L_1^t &
M_2-L_2^t(M_1+K_2)^{-1}L_2}.
\label{s3}
\end{equation}
\end{theorem}

{\sc Proof.} To show that the formula (\ref{kompop})
for the composition of the operators holds,
it is sufficient to verify the identity
$$
\frac1{(2\pi)^{n/2}}
   \int_{\R^n}\exp\Bigl\{ \frac12
   \pmatrix{x^t&y^t} S_1 \pmatrix{x\cr y}\Bigr\}
     \exp\Bigl\{
\frac12 \pmatrix{y^t&z^t} S_2 \pmatrix{y\cr z}\Bigr\}
\,dy
=
$$
\begin{equation}
=
  \lambda(S_1,S_2)
 \exp\Bigl\{ \pmatrix{x^t&z^t} S_3 \pmatrix{x\cr
z}\Bigr\}.
\label{gausss}
\end{equation}
But the integral in the left hand side is
the usual Gauss integral of the form
$$
\int_{\R^n}\exp\bigl\{-\frac12 y^t A y+ y^t b \bigr\}
dy =
 \frac{ (2\pi)^{n/2}}{\det A^{1/2}}
\exp\bigl\{-\frac12 b^t A^{-1}b\bigr\}
$$

{\sc Remark.}
For a symmetric complex $n\times n$ matrix $T$
satisfying
$\Re T>0$, the determinat $\det (T^{-1/2})$
is well defined. Indeed, for each $v\in\C^n$, we have
$\Re\langle T v,v\rangle>0$. Let $v$ be an
eigenvector
of $T$ with the eigenvalue $\lambda$. Then
$0<\Re \langle T v,v\rangle=\Re \lambda
\langle v,v\rangle$
and hence
$\Re \lambda_j>0$. Therefore
we can assume $\det(T^{-1/2})=\prod \lambda_j^{-1/2}$.
\hfill $\square$

\smallskip

Thus we obtain some category $\cK$, whose objects are
0,1,2, \dots, set of morphisms $m\to n$ is
$\Omega_{m,n}$,
and the product of morphisms is given by
 formula (\ref{s3}).

In particular, the  set $\Omega_{n,n}$ is a semigroup.
This semigroup was discussed in the work of Howe
\cite{How}
and the work of Olshanski \cite{Ols} that was
unpublished for
a long time. In particular, it was observed that the
set of
matrices with the invertible block $L$ is
isomorphic to a
subsemigroup
in $\Sp(2n,\C)$ (see Section \ref{remar}).
Nevertheless,
the semigroup $\Omega_{n,n}$ itself
can not be embedded to any group
(since any matrix of the form
$\pmatrix{K&0\cr 0&M}$ is an idempotent).

It appears \cite{NNO}, \cite{Ner1} that the
multiplication (\ref{s3})
hides quite simple algebraic structure,
it is present in next subsection.

\smallskip

{\sc Remark.}
It is natural to consider slightly more
general operators,
whose kernels have the form
\begin{equation}
\exp\bigl\{\frac12\pmatrix{x^t& y^t}) S \pmatrix{x\cr
y}\bigr\}
\delta_L(x,y)
\end{equation}
where $S$ is a {\it nonpositive} definite matrix,
and
$\delta_L$ is a $\delta$-function of some linear
subspace $L\subset \R^{m}\oplus\R^n$.
Obviously, each
kernel of this form is a limit of kernels
of the form (\ref{intop}0.

We restrict ourself to the case described above.
\hfill $\square$

\subsection{Category of Lagrangian linear relations}

Let $V$, $W$ be linear spaces. A {\it linear relation}
$P:V\rr W$ is a linear subspace  in $V\oplus W$.
If $P:V\rr W$, $Q:W\rr Y$ are linear relations,
then their product $QP$ is the linear relation
$V\rr Y$ consisting of vectors
$(v,y)\in V\oplus Y$ satisfying the condition:
 there exists  $w\in W$ such that $(v,w)\in P$,
$(w, y)\in Q$.

{\sc Example.}  Let $A:V\to W$ be a linear
operator.
Then its graph $\graph(A)$ is a linear relation.
If $A:V\to W$, $B:W\to Y$ are linear operators,
then
$$ \qquad\qquad\qquad\qquad
\graph(BA)=\graph(B)\graph(A)
\qquad\qquad\qquad\qquad\qquad\qquad\square
$$

We shall need
the definition of the {\it rank}
of a linear relation
$$\rk P:=\dim P-\dim(P\cap V)-\dim(P\cap W)$$
If $P$ is a graph of a linear operator $A$, then
$\rk P=\rk A$.

\smallskip

 Consider the space
$$V_{n}:=V_{n}^+ \oplus V_{n}^-=\C^n\oplus\C^n$$
We equip this space with two forms,
a skew symmetric bilinear form $L_n$
$$
L_n\bigl((v_+,v_-);(w_+,w_-)\bigr):=
\sum_{j=1}^n (v_+^j w_-^j- v_-^j w_+^j)
$$
 and
a Hermitian form $H_n$
$$
H_n\bigl((v_+,v_-);(w_+,w_-)\bigr):=
\frac 1i
\sum_{j=1}^n (v_+^j \ov w_-^j- v_-^j \ov w_+^j)
$$

{\sc Remark.} Let $g$ be a {\it real} matrix
preserving
the skew-symmetric form
$L_n$.   Obviously, this matrix also preserves
 the form $H_n$. Conversely, let $g$
preserves
the both forms $L_n$, $H_n$. Then it
commutes with complex conjugation and
 hence $g$ is a real matrix. Thus, the group
of linear operators, preserving these two forms
is the real symplectic group $\Sp(2n,\R)$.
\hfill $\square$

\smallskip

{\sc Remark.}
Let us explain how the space $V_n$ with two forms
appears
in a natural way. Consider a space $\R^{2n}$ equiped
with
 a skew-symmetric bilinear form $M$.
Consider its complexification $\C^{2n}$.
We can extend $M$ to $\C^{2n}$ as a bilinear form,
or as a sesquilinear form.
This gives two forms as above.
\hfill $\square$

\smallskip

For each $m$, $n$, we define two forms
on $V_{2m}\oplus V_{2n}$.
$$
L_{m,n}\bigl((v,w),(v',w')\bigr)=
L_m(v,v')-L_n(w,w')
$$
$$
H_{m,n}\bigl((v,w),(v',w')\bigr)=
H_m(v,v')-H_n(w,w')
$$

We define%
\footnote{
As far as the authors know, this category
firstly appeared in \cite{KS}.}
a category $\SP$   whose objects are the
spaces $V_n$.
We define a morphism $V_{n}$  to $V_m$ as a linear
relation $P:V_m\rr V_n$ such that

1. $P$ is Lagrangian%
\footnote{This means that the skew-symmetric form
is zero on $P$ and $P$ has the maximal possible
dimension,
in our case $\dim P=m+n$.}
with respect to the form $L_{m,n}$
In particular $\dim P=m+n$.

2. the form $H_{m,n}$ is   strictly positive definite
on $P$.

We denote the set of all morphisms $m\to n$ by
$\Mor(V_m,V_n)$ or $\Mor_\SP(V_m,V_n)$.

\smallskip

\begin{prop}  Let $P$ be a morphism $V_m$ to
$V_n$,
and $Q$ be a morphism $V_n$ to $V_k$. Then $QP$ is a
morphism
$V_m$ to $V_k$.
\end{prop}

\smallskip

 For proof, see \cite{Ner3},V.1.

In particular, the set $\End(V_n)$
of morphisms  $V_n$ to $V_n$ forms a
semigroup.

As we will observe below in Section \ref{coords}, the
category $\SP$ is
equivalent
to the category $\cK$ defined in the previous
subsection.

\smallskip

{\sc Remark.} Let $g\in\Sp(2n,\R)$.     Then
$\graph(g)$ is Lagrangian with respect to the form
$L_n$
and the form $H_n$ is zero on $\graph(g_n)$.
Thus, the group $\Sp(2n,\R)$ lies on the boundary of
the semigroup $\Mor(V_n,V_n)$.
A variant of the definition of the category $\SP$
including the groups $\Sp(2n,\R)$ is discussed
in \cite{Ner3}.
\hfill$\square$

\subsection{Coordinates on $\Mor(V_m,V_n)$.}
\label{coords}

\begin{lemma}
{\rm a)} {\it Let $P\in\Mor(V_m,V_n)$. Then
there exists an $(n+m)\times(n+m)$ matrix
(V.P.Potapov  transform)
$S=\pmatrix{A & B\cr B^t & C}$
such that  $P$ is the space of solutions
of the following system  of linear equations
$$
\pmatrix{-v^+\cr  w^+}= \pmatrix{A & B\cr B^t & C}
\pmatrix{v^-\cr  w^-}
$$
where
$$v^\pm\in V^\pm_n;\qquad w^\pm\in V^\pm_m$$
Moreover,  $S$ satisfies the conditions

{\rm 1.} $S$ is symmetric ($S=S^t)$

{\rm 2.}   $\Re S<0$}

{\rm b)}
The map $P\mapsto S(P)$ is a bijection
of the set  $\Mor(V_m,V_n)$ to the set of
$(n+m)\times(n+m)$
matrices satisfying  conditions {\rm 1-2}.

{\rm c)}  The product in these coordinates is given by
the
following formula:
if $P\in\Mor(V_m,V_n)$, $Q\in\Mor(V_k,V_m)$, and
$$S(P)=\pmatrix{A & B \cr B^t & C}, \qquad
S(Q)=\pmatrix{K &
L \cr L^t & M}$$
are their Potapov transforms, then
\begin{equation}
S(PQ)=\pmatrix{A-B(C+K)^{-1}B^t & -B(K+C)^{-1}L \cr
-L^t(C+K)^{-1}B^t & M-L^t(C+K)^{-1}L}.
\label{multfor}
\end{equation}
\end{lemma}

{\sc Proof.}
a)- b). The Hermitian form $H_{m,n}$
is zero on $V_n^+\oplus V_m^+$
and it is stricly positive on $P$.
Hence $P\cap (V_n^+\oplus V_m^+)=0$.
Since $\dim P=m+n= \dim  V_n^-\oplus V_m^-$,
the subspace $P$ is the graph of an operator
$ V_n^-\oplus V_m^- \to  V_n^+\oplus V_m^+$.

Our subspace $P$ is Lagrangian with respect to
the skew-symmetric form $L_{m,n}$.
This
is equivalent to     $S=S^t$.

The positivity of the Hermitian form $H_{m,n}$
is equivalent to the condition $\Re S<0$.

\smallskip

c) We have the system of equations
$$
v^+=- Av^- - Bw^-   ; \qquad     w^+= -Kw^-  -Ly^-;
$$
$$
w^+= B^t v^- + Cw^-; \qquad     y^+= L^tw^- + My^-
$$
where
$$
v^\pm\in V_n^\pm;\qquad w^\pm\in V_m^\pm; \qquad
y^\pm\in V_l^\pm
$$
Subtracting the second equation from  the third one,
we obtain
$$
w^-= - (C+K)^{-1} (B^t v_- + L y_-)
$$
Substituting the expression for $w^-$
to the first and the last equations, we get
$$
v^+=-(A-B(C+K)^{-1}B^t) v^-  + B(K+C)^{-1}L    y_-
$$
$$
w^-=- L^t(C+K)^{-1}B^t v^- + (M-L^t(C+K)^{-1}L)y_-.
$$
as it was required.
\hfill$\square$

\subsection {Remark: Linear operators that are
contained
in $\End(V_n)$.}
\label{remar}

 If $n\ne m$, then  $P\in\Mor(V_n,V_m)$
is not a graph of an operator $V_n\to V_m$
(since $\dim P=m+n$).

Let $m=n$. Consider the semigroup $\Gamma_n$
({\it Olshasnki semigroup} \cite{Ols2})
consisting of operators $g:V_n\to V_n$
such that

  a) $g$ preserves the skew-symmetric bilinear form
$L_n$,
i.e., $g\in\Sp(2n,\C)$.

  b) For each nonzero $v\in V_n$,
 $$H_n(gv,gv)<H_n(v,v)$$

If $g\in\Gamma_n$, then $\graph(g) \in\Mor(V_n,V_n)$.

The semigroup $\Gamma_n$ is open (nondense) in
$\Sp(2n,\C)$,
the subgroup $\Sp(2n,\R)\subset\Sp(2n,\C)$
is contained in the closure of $\Gamma_n$.
The semigroup $\Gamma_n$ is open
and dense in $\End(V_n)$.

An element
$P$ of $\End(V_n)$ is contained in $\Gamma_n$
iff $B$ is invertible.

Let $g=\pmatrix{P&Q\cr R&T}\in \Gamma_n$,
i.e.,
$$
\pmatrix{v^+\cr v^-}=
\pmatrix{P&Q\cr R&T}
\pmatrix{w^+\cr w^-}
$$
Expressing $v^+$, $w^+$ in terms of
$v^-$, $w^-$, we obtain that
the Potapov transform of $g$  is
$$S=
\pmatrix{-PR^{-1}& -Q+PR^{-1}T\cr
         R^{-1}& -R^{-1}T}
$$

\subsection{Explicit correspondence between Lagrangian
linear
relations and Gauss operators.}
Let $P\in\Mor(V_n,V_m)$. Let
$$S(P)=\pmatrix{K&L\cr L^t&M}$$
be its Potapov transform.
Consider the integral operator
$$\B[S(P)]:L^2(\R^n)\to L^2(\R^m)$$
given by
$$
(\B [S(P)]f)(x)=\frac
1{(2\pi)^{(n+m)/4}}\int_{\R^n}\exp\{
\frac{1}{2}\pmatrix{x^t
& y^t}S(P)\pmatrix{x \cr y}\}
f(y)dy.
$$

\begin{theorem}  Let $P\in\Mor(V_n,V_m)$,
$Q\in\Mor(V_m,V_l)$, and
$S(P)=\pmatrix{A & B \cr B^t & C}$, $S(Q)=\pmatrix{K &
L \cr L^t & M}$ be their Potapov transforms.
Then
$$
\B[S(Q)]\B[S(P)]=\lambda(Q,P) \B[S(QP)]
$$
where

$$
\lambda(P,Q)=\det(-C-K)^{-1/2}.
$$
\end{theorem}

This theorem is a corollary of the multiplication
formula
(\ref{multfor})
and  formulae (\ref{kompop})-(\ref{s3})
 for the product of integral operators.

{\sc Remark} We can say that $P\mapsto \B[S(P)]$ is a
projective
representation
of the category $\SP$, on definitions of
representations of
categories see \cite{Ner3}.

\smallskip

{\sc Remark.}
As we mention above, in \cite{Ner3}
another definition
of category $\SP$, which gives a slightly larger
structure, is used.
This category is equivalent to enlarged category of
Gauss operators mentioned in 2.1. Explicit formulae
for the correspondence between elements of
 $\End(V_n)$ and Gauss operators are contained in
\cite{Hor}.

\subsection{Heisenberg algebra and another description
of the correspondence between linear relations and
Gauss
operators}
\label{heisalg}

To each
$$
\alpha=(\alpha^+;\alpha^-)^t=(\alpha^+_1,\dots,
\alpha^+_n;
\alpha^-_1,\dots, \alpha ^-_n)^t\in V_n=V_n^+\oplus
V_n^-
$$
we associate an operator
$$
\widehat A(\alpha)=
\sum \alpha^+_j x_j+\sum  \alpha^-_j \frac\partial
{\partial x_j}.
$$
All operators of this type form the
complex Heisenberg
algebra.

\begin{prop}
Let $P\in \Mor(V_n,V_m)$.

a) For each $(\alpha,\beta)\in P$,
$$
\widehat A(\alpha) \B[S(P)]= \B[S(P)]\widehat A(\beta)
$$

b) Let $R:\cD(\R^n)\to
\cD'(\R^m)$ be a bounded operator
satisfying the equality
$$
\widehat A(\alpha) R= R\widehat A(\beta)
$$
for all $(\alpha,\beta)\in P$.
Then $R=\lambda \B [S(P)]$ for some $\lambda\in\C$.
\label{operat}
\end{prop}

\smallskip

{\sc Proof.}
a) It is sufficient to prove that the kernel
$K=K(x,y)$ of $\B [P]$
satisfies
$$
\Bigl(\sum \alpha^+_j x_j+\sum  \alpha^-_j
\frac\partial
{\partial x_j}         \Bigr)
 \int_{\R^n} K(x,y) \,f(y)\,dy=
$$
$$
= \int_{\R^n} K(x,y)\,
\Bigl(\sum \beta^+_k y_k+\sum  \beta^-_k \frac\partial
{\partial y_k}    \Bigr)
 \,f(y)\,dy
$$

This is equivalent to
 the following system of differential
equations:
$$
\sum_{j=1}^m(\alpha_j^+x_j+\alpha_j^-\frac{\partial}{\partial
x_j}) \exp\Bigl\{ \frac12 \pmatrix {x^t&y^t}
\pmatrix{A&B\cr B^t &C} \pmatrix {x \cr y}\Bigr\} =
$$
$$
=
\sum_{k=1}^n(\beta_k^+ y_k -
\beta_k^-\frac{\d}{\d y_k}\Bigl)
\exp\Bigl\{ \frac12 \pmatrix {x^t&y^t}
\pmatrix{A&B\cr B^t &C} \pmatrix {x \cr y} \Bigr\}
$$
After differentiation, we obtain
\begin{equation}
\alpha^+=-A\alpha^- -B\beta^-;\qquad
\beta^+=B^t\alpha^-+C\beta^-
\label{xxx}
\end{equation}

b) By a Kernel Theorem   (see \cite{Hor0}, 5.2.6)
the operator $R$ has the form
$$
Rf(x)=\int_{\R^n} L(x,y)\,f(y)\, dy
$$
where $L\in \cD'(\R^n\times\R^n)$.
The distribution $L$ must satisfy
\begin{equation}
\label{eqL}
\sum_{j=1}^m(\alpha_j^+x_j+\alpha_j^-\frac{\partial}{\partial
x_j}) L(x,y)              =
\sum_{k=1}^n(\beta_k^+y_k-\beta_k^-\frac{\d}{\d
y_k}\Bigl)
L(x,y)
\end{equation}
for all $\alpha$, $\beta$.

Consider the distribution
$$
M(x,y):=L(x,y)
\exp\Bigl\{-\frac12 \pmatrix {x^t&y^t}
\pmatrix{A&B\cr B^t &C} \pmatrix {x \cr y} \Bigr\}
$$
The system (\ref{eqL})
implies
$$
\sum_{j=1}^m(\alpha_j^-\frac{\partial}{\partial
x_j}) M(x,y)
=
\sum_{k=1}^n(\beta_k^-\frac{\d}{\d
y_k}\Bigl)
M(x,y)
$$
for all $\alpha_j^-$, $\beta_k^-$.
Hence $M(x,y)$ is a constant.
\hfill $\square$

\smallskip

{\sc Remark.} Let $n=0$. Let $P\in\Mor_\SP(V_0,V_n)$.
Its Potapov transform is $(m+0)\times(m+0)$ matrix
$A$. The corresponding Gauss operator
$$
{\cal B}(P):L^2(\R^0)=\R\to L^2(\R^m)
$$
is given by
$$
s\mapsto s\cdot \exp\bigl\{\frac12 x^t A x\bigr\}
$$

\section{Integral operators with theta-kernels}
\label{opertheta}

\setcounter{equation}{0}

\subsection{Zak transform and its properties}

Consider  tori $\T^{n}=\R^{n}/(2\pi \Z)^{n}$,
$\T^{2n}=\R^{2n}/(2\pi \Z)^{2n}$.

Denote by $G(\R^{2n})$  the subspace of
$C^{\infty}(\R^{2n})$ that
consists of functions $g(x,\xi)$
with the properties
\begin{equation}
g(x,\xi+2\pi k)=g(x,\xi),\qquad g(x+2\pi
k,\xi)=e^{-ik\cdot \xi}
g(x,\xi)
\label{transflaw}
\end{equation}
 where $x\in\R^n$, $\xi\in \R^n$, $k\in\Z
^n$. We equip this space
 with the topology of uniform convergence with
all derivatives on $[0,2\pi]^n$.

We also define an inner product
in $G(\R^{2n})$
by the formula
$$
\langle g_1,g_2\rangle= \frac{1}{(2\pi )^n}
\int_{[0,2\pi]^{2n}} g_1(x,\xi)\ov{g_2(x,\xi)} dx d\xi
$$

{\sc Remark.}
In this formula we can replace integration
over the cube $[0,2\pi]^{2n}$
by  integration
over an arbitrary fundamental domain of the lattice
$\Z^{2n}$ in
$\R^{2n}$.
\hfill$\square$

The completion of $G(\R^{2n})$ with respect
to this inner product is
identified with $L^2$ on the cube $[0,2\pi]^{2n}$.

\smallskip

{\sc Remark.} Obviously, we can naturally identify the
space $G(\R^{2n})$
with the space of smooth sections of
a certain line bundle
$\cL \to \T^{2n}$. Indeed, $\cL=\R^{2n}\times \C /
\sim$, where the equivalence
relation is $(x,\xi,\zeta)\sim (x+2\pi k,\xi+2\pi m,
e^{ik\cdot\xi}\zeta)$
for any $x,\xi\in\R^n$, $m,k\in\Z^n$, $\zeta\in\C$.
\hfill$\square$

For a function $f$ on $\R^n$ we define its {\it Zak
transform}%
\footnote{Another term is also used: {\it Weil--Brezin
transform}, see \cite{Bre}.}
$\fZ_n f$ by the formula
$$
\fZ_n:f(x)\mapsto g(x,\xi)=\sum_{k \in \Z^n}f(x+2\pi
k)e^{ik\cdot \xi}.
$$
where $\xi\in\R^n$.
It is easy to see that the function $g(x,\xi)$
satisfies
(\ref{transflaw}).

\begin{theorem}[see \cite{KG}, Problems 475, 666,
\cite{Fol}]
{\rm a)}  The
Zak transform is a bounded invertible operator
\begin{equation}
\fZ_n:\S (\R^{n})\to G(\R^{2n})
\label{mapfi}
\end{equation}

{\rm b)}  The inverse transform is given by
$$
(\fZ_n^{-1}g)(x)=\frac{1}{(2\pi)^n}\int_{\T^n}g(x,\xi)d\xi.
$$

{\rm c)}  The Zak transform is a unitary operator
$L^2(\R^{n})\to L^2([0,2\pi]^{2n})$.
\end{theorem}

\smallskip

{\sc Remark.}
The Zak transform has the following property
unusual for the classical theory
of integral transforms: it changes a functional
dimension; i.e., it identifies a space of functions
of $n$ variables and a space of functions of $2n$
variables.

\smallskip

\begin{prop}
[\cite{Car}]
The Zak transform
maps the operator $f\mapsto x_jf$ in $L^2(\R^n)$
to the operator $g(x,\xi)\mapsto
(x_j+\frac{2\pi}{i}\frac{\d}{\d\xi_j})$
and the operator $\partial/\partial x_j$
 to $\partial/\partial x_j$.
\label{troper}
\end{prop}

{\sc Remark.}
In other words, the Zak transform intertwines the
standard representation
of the Heisenberg algebra in $\cS(\R)$ given by
the operators $x_j$, $\partial/\partial x_j$
and
the representation in
$G(\R^{2n})$ given  by the operators
$\frac{1}{i}\frac{\d}{\d x_j}$ and
$x_j+\frac{2\pi}{i}\frac{\d}{\d\xi_j}$, $j=1,...n$,
see also \cite{LV}, \cite{Mum3}, Section 2.
\hfill $\square$

\subsection{$\theta$-functions,} see also
\cite{Mum3}, Section 2.
For a negative definite $a\times a$
 symmetric matrix $T$ define a
theta-function
$$
\theta [T;z,\zeta]=\sum_{k\in\Z^a}\exp \{\frac12
(z+2\pi
k)^t T
(z+2\pi k)  \}e^{ik\cdot \zeta},
$$
$z,\zeta\in \R^a$,
clearly $\theta [T;z,\zeta]$ is the image $\fZ_a f$ of
the Gaussian $f(z)=\exp\{\frac12 z^tTz\}$
under the Zak transform.

By unitarity of the Zak transform,
\begin{equation}
\langle \theta [T_1;z,\zeta],\theta [T_2;z,\zeta]
\rangle=
\int_{\R^a}\exp\{ \frac{1}{2} z^t(T_1+T_2)z  \}dz=
\frac{ (2\pi)^{\frac{a}{2}}}
{\det (-T_1-T_2)^{-1/2}}.
\end{equation}
(below we use a modified variant of this identity
in  Section \ref{proof1}).

The Gaussian $\exp\{\frac12 x^tTx\}$
satisfies the differential equation
$$
\sum_{j=1}^a \bigl(\gamma^+_j z_j+
   \gamma^-_j\frac{\partial} {\partial z_j}\bigr)
     \exp\{\frac12 z^tTz\}=0
\qquad \mbox{if} \quad\gamma^+=-T \gamma^- .
$$
By Proposition \ref{troper} the corresponding
$\theta$-function
satisfies the equation

\begin{equation}
\sum_{j=1}^a \bigl(\gamma^+_j
(z_j+\frac{2\pi}{i}\frac{\d}{\d\zeta_j})+
   \gamma^-_j\frac{\partial} {\partial z_j}\bigr)
     \theta [T;z,\zeta]=0.
\label{difurtheta}
\end{equation}

\subsection{Formulation of result}

For a symmetric $(m+n)\times (m+n)$ matrix $S$ denote
$b(x,y)=\frac{1}{2}\pmatrix{x^t & y^t}S\pmatrix{x \cr
y}$, $x\in\R^m$, $y\in\R^n$,
define a theta-kernel $\K_{n,m}[S](x,\xi;y,\eta)$
as
$$
\K_{n,m}[S](x,\xi;y,\eta)=\sum_{k\in\Z^m}\sum_{l\in\Z^n}
\exp\{ b(x+2\pi k,y +2\pi l)\} e^{ik\cdot
\xi}e^{-il\cdot
\eta},
$$
and the corresponding integral operator
$$\cQ_{n,m}[S]:\,G(\R^{2n})\to
G(\R^{2m})$$
by
$$
(\cQ_{n,m}[S]g)(x,\xi)=
\frac{1}{(2\pi)^{3(n+m)/4}}
\int_{\T^{2n}}\K_{n,m}[S](x,\xi;y,\eta)g(y,\eta)dyd\eta.
$$
where $g\in G(\R^{2n})$.

\begin{theorem}
\label{main1}
Let $P_1\in Mor_\SP(V_{q},V_m)$, $P_2\in
Mor_\SP(V_{n},V_{q})$ be linear relations.
Let
\begin{equation}
S(P_1)=\pmatrix{A & B \cr B^t & C},
\qquad
S(P_2)=\pmatrix{K & L \cr L^t & M}
\end{equation}
be their Potapov transforms.
Then the following formula holds for the composition
\begin{equation}
\cQ_{q,m}[S(P_1)] \cQ_{n,q}[S(P_2)]=
\det(-C-K)^{-1/2} \cQ_{n,m}[S(P_1 P_2)].
\label{lambq}
\end{equation}

\smallskip

In particular, we obtain a projective representation
of the
category {\bf Sp}.
\end{theorem}
\smallskip

\begin{theorem}
$
\fZ_m B[S] \fZ_n^{-1}=\cQ_{n,m}[S]
$
\label{main2}
\end{theorem}
\smallskip

The rest of the paper is the proof of these theorems.

\subsection{Computation of the kernel. Proof of
Theorem \ref{main2}}
\label{compker}

Now we shall prove that
the operator
$$ Q:\fZ _m \B[S]\fZ_n^{-1}:G(\R^{2n})\to G(\R^{2m})$$
is equal {\it up to a scalar factor} to the operator
with the kernel
$\frac{1}{(2\pi)^{3(n+m)/4}}\K_{n,m}[S]$, 
where
$$
\K_{n,m}[S](x,\xi;y,\eta)=
\sum_{k\in\Z^m}\sum_{l\in\Z^n}
K(x+2\pi k,y +2\pi l)e^{ik\cdot \xi}e^{-il\cdot\eta},
$$
where $K(x,y)=\exp \{ \frac{1}{2}
\pmatrix{x^t & y^t}S\pmatrix{x \cr y} \}$,
i.e. $\frac{1}{(2\pi)^{(n+m)/4}}K(x,y)$ is the
kernel of $\B [S]$. In Section \ref{endofproof} we
will show
that the scalar factor is 1.

First, $\K$ satisfies the quasiperiodicity
conditions
\begin{eqnarray*}
{\K}_{n,m}[S](x+2\pi k,\xi;y+2\pi l,\eta)&=&
e^{-i \xi\cdot k+i\eta\cdot
l}{\K}_{n,m}[S](x,\xi;y,\eta),
\\
{\K}_{n,m}[S] (x,\xi+2\pi k;y,\eta+2\pi l)&=&
{\K}_{n,m}[S](x,\xi;y,\eta).
\end{eqnarray*}
and hence $\K$ is the kernel of an operator
$G(\R^n)\to G(\R^m)$.

Recall that for $(\alpha,\beta)\in P$, or
equivalently,
for $\alpha$, $\beta$ satisfying (\ref{xxx})
\begin{equation}
\hat A(\alpha)\B[S]=\B[S]\hat A(\beta),
\label{intertw}
\end{equation}
where
$$
\hat
A(\alpha)=\sum_{j=1}^m(\alpha_j^+x_j+\alpha_j^-\frac{\d}{\d
x_j}),
\qquad
\hat
A(\beta)=\sum_{j=1}^n(\beta_j^+y_j+\beta_j^-\frac{\d}{\d
y_j}),
$$
and we showed
that $K(x,y)$ satisfies the differential
equation
\begin{equation}
\sum_{j=1}^m(\alpha_j^+x_j+\alpha_j^-\frac{\d}{\d
x_j})K(x,y)=
\sum_{j=1}^n(\beta_j^+y_j-\beta_j^-\frac{\d}{\d
y_j})K(x,y).
\label{eqK}
\end{equation}

Applying the Zak transform to (\ref{intertw}) we get:
$$
(\fZ_m \hat A(\alpha)\fZ_m^{-1}Qg)(x,\xi)=(Q\fZ_n \hat
A(\beta)\fZ_n^{-1}g)(x,\xi),
$$
where $g\in G(\R^{2n})$.
Hence
the kernel $\K$ of the operator $Q$ satisfies the
identity
\begin{eqnarray*}
\fZ_m \hat A(\alpha)\fZ_m^{-1}\int_{[0,2\pi]^{2n}}
\tilde{\K}(x,\xi;y,\eta)g(y,\eta)dy d\eta=
\\
\int_{[0,2\pi]^{2n}}\tilde{\K}(x,\xi;y,\eta)
\fZ_n \hat A(\beta)\fZ_n^{-1}g(y,\eta)dy d\eta .
\end{eqnarray*}
By Proposition   \ref{troper}
\begin{eqnarray*}
\fZ_m \hat A(\alpha)\fZ_m^{-1}=
\sum_{j=1}^m
(\alpha_j^+(x_j+\frac{2\pi}{i}\frac{\d}{\d\xi_j})+\alpha_j^-\frac{\d}{\d
x_j}),
\\
\fZ_n \hat A(\beta)\fZ_n^{-1}=
\sum_{j=1}^n
(\beta_j^+(y_j+\frac{2\pi}{i}\frac{\d}{\d\eta
_j})+\beta_j^-\frac{\d}{\d y_j}).
\end{eqnarray*}
Hence
$\tilde{\K}(x,\xi;y,\eta)$ satisfies the
differential
equations
\begin{eqnarray*}
\sum_{j=1}^m
(\alpha_j^+(x_j+\frac{2\pi}{i}\frac{\d}{\d\xi_j})+\alpha_j^-\frac{\d}{\d
x_j})
\tilde{\K}(x,\xi;y,\eta)=
\\
=\sum_{j=1}^n
(\beta_j^+(y_j-\frac{2\pi}{i}\frac{\d}{\d\eta
_j})-\beta_j^-\frac{\d}{\d y_j})
\tilde{\K}(x,\xi;y,\eta)
\end{eqnarray*}
and if $\tilde{\K}={\rm const}\cdot\K_{n,m}[S]$ then
they are satisfied because of (\ref{difurtheta})
with $a=n+m$,
$$\gamma^+=\pmatrix{\alpha^+ \cr -\beta^+},\quad
\gamma^-=\pmatrix{\alpha^- \cr \beta^-},\quad
z=\pmatrix{x \cr y},\quad \zeta=\pmatrix{\xi \cr
-\eta},\quad T=S.$$
Note that the condition $\gamma^+=-T\gamma^-$
is satisfied by (\ref{xxx}).

Finally we note that the operator $\B[S]:L^2(\R^n)\to
L^2(\R^m)$ satisfying (\ref{intertw}) is unique
up to a constant factor
by Proposition \ref{operat} (b),
therefore, clearly,  the corresponding operator
$Q:G(\R^{2n})\to G(\R^{2m})$ is unique up to a scalar
factor
too.

\subsection{Composition formula. Proof of Theorem
\ref{main1}.}
\label{proof1}

Let $Q_1:=\cQ_{q,m}[S_1]$ and $Q_2:=\cQ_{n,q}[S_2]$ be
two
operators
with the kernels $\frac{1}{(2\pi)^{3(q+m)/4}}\K_1$ and
$\frac{1}{(2\pi)^{3(n+q)/4}}\K_2$ respectively, where
$$
\K_1=\K_1(x,\xi;y,\eta)
=\K_{q,m}[S_1](x,\xi;y,\eta)=
=\sum_{k\in\Z^m}\sum_{l\in\Z^{q}}
e^{b_1(x+2\pi k,y +2\pi l)}e^{ik\cdot \xi}e^{-il\cdot
\eta},
$$
$$
\K_2=\K_2(x,\xi;y,\eta)=
\K_{n,q}[S_2](x,\xi;y,\eta)=
\sum_{k\in\Z^{q}}\sum_{l\in\Z^n}
e^{b_2(x+2\pi k,y +2\pi l)}e^{ik\cdot \xi}e^{-il\cdot
\eta},
$$
and $b_1(.,.)$ and $b_2(.,.)$ are the quadratic
forms
associated to $S_1$, $S_2$.

Our Theorem \ref{main1} is equivalent to

\begin{prop}
\label{dubl}
 The composition $Q_3=Q_1Q_2$
is the
operator
defined by
$$
(Q_3g)(x,\xi)=
\frac{1}{(2\pi)^{3(n+m)/4}}
\int_{[0,2\pi]^{2n}} \K
_3(x,\xi;y,\eta)g(y,\eta)dyd\eta,
$$
where $\K _3(x,\xi;y,\eta)$ is  given by
\begin{equation}
\!\!\!\!\!\!
\lambda (S_1,S_2)
\sum_{k\in\Z^m}\sum_{l\in\Z^n}
\exp \bigl\{\frac{1}{2} \pmatrix{x^t+2\pi k^t &
y^t+2\pi
l^t}
S_3\pmatrix{ x+2\pi k \cr y+2\pi l}\Bigr\}
e^{ik\cdot \xi}e^{-il\cdot \eta},
\label{kerkomp}
\end{equation}
where $\lambda(S_1,S_2)$ is given by (\ref{lamb})
and $S_3$ is given by (\ref{s3}).
\end{prop}

{\bf Proof.}
We have:
$$
(Q_3g)(x,\xi)=(Q_1Q_2g)(x,\xi)=
$$
$$
\frac{1}{(2\pi)^{3(n+m+2q)/4}}
\int_{[0,2\pi]^{2q}}\int_{[0,2\pi]^{2n}}
\K_1(x,\xi;s,\zeta)\K_2(s,\zeta;y,\eta)g(y,\eta)dyd\eta
dsd\zeta,
$$
hence
$$
(2\pi)^{3q/2}\K_3(x,\xi;y,\eta)=
\int_{[0,2\pi]^{2q}}
\K_1(x,\xi;s,\zeta)\K_2(s,\zeta;y,\eta) dsd\zeta=
$$
$$
\int\limits_{[0,2\pi]^{2q}}
\sum_{k\in\Z^m}\sum_{p\in\Z^{q}}
e^{b_1(x+2\pi k,s +2\pi p)}e^{ik\cdot \xi}e^{-ip\cdot
\zeta}
\sum_{r\in\Z^{q}}\sum_{l\in\Z^n}
e^{b_2(s+2\pi r,y +2\pi l)}e^{ir\cdot
\zeta}e^{-il\cdot \eta}
dsd\zeta=
$$
$$
\sum_k\sum_l e^{ik\cdot \xi}e^{-il\cdot
\eta}\int_{[0,2\pi]^{2q}}
\sum_p\sum_r e^{b_1(x+2\pi k,s +2\pi p)}
e^{b_2(s+2\pi r,y +2\pi l)}e^{ir\cdot
\zeta}e^{-ip\cdot\zeta}dsd\zeta .
$$
For fixed $x$, $y$, $\xi$, $\eta$
the integral above is $(2\pi)^q$ times the inner
product (in
$G(\R^{2q})$)
 of the Zak transforms of two Gaussian functions,
therefore
 the expression above becomes
$$
\frac{1}{(2\pi)^{q/2}}
\sum_k\sum_l e^{ik\cdot \xi}e^{-il\cdot
\eta}
\Bigl\langle \fZ_{q} (e^{b_1(x+2\pi k,s)}),\fZ_{q}
(e^{b_2(s,y +2\pi l)})
\Bigr\rangle_{L^2([0,2\pi]^{2q})} .
$$
The Zak transform (\ref{mapfi}) is unitary.
Hence we can replace the inner product in
$L^2[0,2\pi]^{2q}$
by the inner product in $L^2(\R^{q})$.
Therefore  we obtain:
$$
\K_3(x,\xi;y,\eta)=
\frac{1}{(2\pi)^{q/2}}
\sum_k\sum_l e^{ik\cdot
\xi}e^{-il\cdot \eta}
\int_{\R^{q}}e^{b_1(x+2\pi k,s)} e^{b_2(s,y+2\pi
l)}ds.
$$
The last integral is the integral
(\ref{gausss}) with $x,y,z,n$ replaced by
$x+2\pi k$, $s$, $y+2\pi l$, $q$,
and this finishes the proof.

\subsection{End of proof of Theorem \ref{main1}}
\label{endofproof}

In section \ref{compker} we evaluated the kernel of
the operator
$
\fZ_m \B[S] \fZ_n^{-1}
$
up to a scalar factor, i.e.,
\begin{equation}
\fZ_m \B[S] \fZ_n^{-1}=\sigma(m,n;S) \cQ_{n,m}[S]
\label{sigma1}
\end{equation}
where  $\sigma(m,n;S)\in \C$.
We intend to prove that
\begin{equation}
\sigma(m,n;S)=1 \qquad \mbox{\it for all} \quad m,\,\,
n,\,\, S.
\label{1}
\end{equation}

A priori we know the following facts
about the function  $\sigma$.

\begin{lemma}
{\rm a)}
$\sigma(m,n;S)$  is a holomorphic function in
the variable
$S$.

$$\!\!\!\!\!\!\!\!  \!\!\!\!\!\!\!\! \!\!\!\!\!\!\!\!
\!\!\!\!\!\!\!\!
b) \qquad \sigma(0,m;S)=1
\qquad\qquad\qquad\qquad\qquad\qquad\qquad
;$$

\begin{equation}
\!\!\!\!\!\!\!\!\!\!\!\!\!\!\!\!\!\!\!\!\!\!\!\!\!\!\!\!\!\!\!\!
{\rm c)}  \qquad\qquad
\sigma(m,n;S(P_1) \sigma(k,m;S(P_2)  =
\sigma(k,n;S(P_1P_2))
\label{partc}
\end{equation}    \qquad\qquad\qquad\qquad\qquad\qquad
\end{lemma}

{\sc Proof.}
a) Indeed,  the both part of the equality
(\ref{sigma1}) are
holomorphic in
$S$.

b) Indeed, the Gauss operators corresponding to
elements $\Mor_\SP(V_0,V_m)$ are described at the end
of Section \ref{heisalg}.
Due to Section 3.2, we know their images with respect
 to the Zak transform
exactly, not up to a constant multiplicative factor.

c)
We observe (see (\ref{lamb}), (\ref{lambq})) that
in the identities
$$
\cB[S(P_1)] \cB[S(P_2)]=\lambda(P_1,P_2)  \cB[S(P_1
P_2)]
$$
$$
\cQ_{m,n}[S(P_1)]  \cQ_{k,m}(S(P_2))=
\lambda(P_1,P_2)
 \cQ_{k,n}(S(P_1 P_2))
$$
the scalar factors  $\lambda(P_1,P_2)$
coincide.  This proves (\ref{partc}).
\hfill $\square$

\smallskip

Below we easily  reduce (\ref{1}) to our lemma.

{\it Step 1.}
For each idempotent $P\in\End(V_n)$,
the identity (\ref{sigma1}) implies
$\sigma(n,n;S(P))=1$.
Thus we need the description of idempotents
in $\End_\SP(V_n)$.

\begin{lemma}
{\rm a)}  An element $P\in\End_{\SP}(V_n)$
is an idempotent (i.e. $PP=P$) iff $P\subset V_n\oplus
V_n$
has the form
$P=Y\oplus Z$, where $Y$ is a linear subspace in the
first copy
of $V_n$ and $Z$ is a linear subspace in the second
copy of
$V_n$.

{\rm b)}  Let $P\in \End(V_n)$ be an idempotent. Then
for
each $Q\in \End(V_n)$, the linear relation $QP$ is an
idempotent.
\end{lemma}

\smallskip

The proof of a) is
a straightforward verification.
The statement b) is a corollary of a).

{\it Step 2.}
Considering an idempotent $
P\in\End(V_n)$
and arbitary $Q\in\End(V_n)$, we obtain
$$
\sigma(n,n);S(Q))=\sigma(n,n;S(QP))/\sigma(n,n;S(P))=1/1
$$
i.e.,
\begin{equation}
 \sigma(n,n;S(Q))=1 \qquad \mbox{for all
$Q\in\End(V_n).$}
\label{sigmnn}
\end{equation}

{\it Step 3.}
Let $n>m$.
Now fix $R\in\Mor(V_n,V_m)$ of the maximal possible
rank.
Let $P$ range in $\End(V_n)$. Then the set of all
possible
products $RP$ is open (nondense) in
$\Mor(V_n,V_m)$.

By (\ref{partc}),(\ref{sigmnn})
we have:
$$
\sigma(n,m;S(RP))=
\sigma(n,n,P)\sigma(n,m;S(R)) =\sigma(n,m;S(R)).
$$
Hence the expression  $\sigma(n,m;S(T)) $ is a
constant on
the open
subset consisting of products $RP$.
But $\sigma(n,m;S(R))$ is holomorphic and hence
it is some constant $\sigma_{n,m}$.

For $m>n$ we repeat the same considerations with
products
$PR$ with fixed $R\in\Mor(V_n,V_m)$ and $P$ ranging in
$\End(V_n)$.

{\it Step 4.}
Now we have
$$
\sigma_{0,n}\sigma_{n,m}=
\sigma_{0,m}.
$$
But $\sigma_{0,k}=1$ for all $k$.
Thus $\sigma_{m,n}=1$.
This finishes proof of the equality (\ref{1}).

{

\end{document}